\documentclass[]{article}

\usepackage{makeidx}         %
\usepackage{graphicx}        %
            
\usepackage{multicol}        %
\usepackage[bottom]{footmisc}%

\usepackage{newtxtext}       %
\usepackage[varvw]{newtxmath}       %
\usepackage{authblk}
\usepackage{svg}

\usepackage{hyperref}

\usepackage{tikz}

\newcommand{\define}[1]{{\bf \boldmath{#1}}}
\newcommand{\VR}{\mathrm{VR}}
\newcommand{\length}{\mathrm{length}}
\newcommand{\SP}{\mathrm{d}_{\mathrm{SP}}}

\newtheorem{definition}{Definition}

\newtheorem{example}{Example}

\makeindex             %

\date{}

\begin{document}

\title{Fractal dimensions of complex networks: advocating for a topological approach}

\author[1]{Rayna Andreeva}
\author[2]{Hayde\'e Contreras-Peruyero}
\author[3]{Sanjukta Krishnagopal}
\author[4]{Nina Otter}
\author[5]{Maria Antonietta Pascali}
\author[6]{Elizabeth Thompson}

\affil[1]{University of Edinburgh, School of Informatics, United Kingdom}
\affil[2]{ Centro de Ciencias Matem\'aticas, UNAM, M\'exico}
\affil[3]{University of California Santa Barbara, USA}
\affil[4]{DataShape, Inria-Saclay,  France}
\affil[5]{Institute of Information Science and Technologies of the National Research Council, Italy}
\affil[6]{Washington State University, USA }

\maketitle

\abstract{Topological Data Analysis (TDA) uses insights from topology to create representations of data able to capture global and local geometric and topological properties. Its methods have successfully been used to develop estimations of fractal dimensions for metric spaces that have been shown to outperform existing techniques. In a parallel line of work, networks are ubiquitously used to model a variety of complex systems. Higher-order interactions, i.e., simultaneous interactions  between more than two nodes, are wide-spread in social and biological systems, and simplicial complexes, used in TDA,  can capture important structural and topological properties of networks by modelling such higher-order interactions. In this position paper, we advocate for  methods from TDA to be used to estimate fractal dimensions of complex networks, we discuss the possible advantages of such an approach and outline some of the challenges to be addressed. }

\section{Networks as relational models for data}

Many real-world phenomena, such as social systems \cite{holland1976local, brandes2013social}, biological processes \cite{pavlopoulos2011using,bullmore2009complex} or  communication systems \cite{hajij2018visual}, amongst others, are characterised by being constituted of entities and %
pairwise relationships between them. 
Network science uses graphs to model such systems and studies properties of the resulting graphs, usually called \emph{complex networks} or simply \emph{networks}.\footnote{The difference between the terms ``graph'' and ``network'' is thus in the eyes of the beholder: while both terms refer to the same underlying mathematical object, the problems and questions studied in the fields of  network science and graph theory are very different. While in graph theory one studies properties of graphs considered as abstract mathematical objects, in network science one is instead interested in studying properties of the subset of graphs that arise as models of real-world systems. }

Complex networks have been studied since at least the 1930s. One of the first applications of networks, which has driven much of the early development of the field, was developed by structural social scientists \cite{freeman2004development}, who were motivated by the desire to model social systems not by studying social actors --- which was and still is the common approach in many branches of social science research --- , but rather  the relationships between them.
More generally, weighted networks provide a suitable model for systems with interacting agents where interactions have different strengths or flavours, and edge-weights capture, e.g., the strength of interaction between nodes.
The field of network science received growing interest from researchers outside the field from the early 2000s, when it was demonstrated that  networks can 
provide a suitable model for the study of dynamical processes such as the spread of a disease in a population \cite{PSV01}. The field continues to evolve and finds applications in a variety of domains, from transportation networks \cite{Pitts-1995}, to growth of fungi \cite{dikec2020hyphal} or  modelling the spread of the COVID-19 disease \cite{lee2024uncovering}, to cite a few applications. 
We refer to \cite{newman} for an overview of the main methods and applications in the field.

\section{What is the dimension of a network?}

The dimension of a process, an object, or a shape, can roughly be thought of as the minimal number of parameters, or coordinates, needed to describe it. While dimension is a well-established notion in  geometry and topology,  the question of what the dimension of a network is has been addressed only relatively recently.

A first approach one could take in defining the dimension of a network is to embed it into a space with a well-known notion of dimension, such as Euclidean or hyperbolic space. 
A problem with this approach is that for an arbitrary network there is in general not a unique choice of embedding into Euclidean or hyperbolic space, and thus different choices may result in different values for the dimension. The study of network embeddings is a very active field of research, and we point the reader to \cite[Chapter 17.7.4]{newman} and references therein for further reading, and to \cite{DKKPT22} for a recent study comparing  network embeddings. 

The question of what an \emph{intrinsic} dimension of a network is was, to the best of our knowledge, first addressed in the 1990s. Motivated by the hypothesis that physics, on the  Planck scale, has to be described through discrete concepts, Requardt introduced  a collection of networks to model physical vacuum and space-time  \cite{Manfred_Requardt_1998}. The desire to understand ``what properties actually are encoded in a notion like dimension on the most fundamental physical level'' \cite{NR1998} then lead to the 
need  to define   a  notion of dimension that takes into account the \emph{discrete} nature of such networks,  and does not rely on the dimension of a \emph{continuous} space in which they  might be embedded, see \cite{NR1998} and  references therein. 
In the following decades, 
 many different methods to define intrinsic dimensions of networks  have been developed;
most of these methods are closely related to known  so-called \emph{fractal dimensions} from geometry and topology, such as the box-counting dimension and Hausdorff dimension, which generalise  usual topological notions of dimension to beyond the integers. 
The study of an intrinsic notion of dimension inscribes itself in the more general study of intrinsic properties characterising complex networks, such as diameter, or the number of connected components.

\subsection{Box-counting dimension of subsets of Euclidean space}\label{SS:box-counting}
The \emph{box-counting dimension} (also known as \emph{Minkowski dimension}) relies on covering a space with covering sets of a specific size, and in studying how the minimal number of sets  needed to cover the space varies as one changes the size parameter of the sets. More precisely, one gives the following definition:

\begin{definition}
Let $S\subset \mathbb{R}^n$ be a non-empty and bounded subset. Let $N(\epsilon)$ be the minimal  number of closed balls of radius $\epsilon$ that cover $S$. The \define{box-counting dimension}\footnote{The name ``box-counting'' is evocative of an equivalent definition that one may give, by considering cubes (i.e., ``boxes'') of fixed side length. There exist at least five different types of covering sets that give equivalent definitions of box-counting dimension, see \cite[Section 3.1]{falconer2004fractal}.} of $S$ is defined as the following limit, if it exists:
\[
\dim_{box}(S):=\lim_{\epsilon\to 0} \frac{\log( N(\epsilon))}{\log(1/\epsilon)}\, .
\]

\end{definition}

Given   an interval in $\mathbb{R}$, if we halve the radius of the balls needed to cover it, then we will need twice as many balls. Thus, asymptotically, for the radius of the balls going to zero, we have that the minimal number of covering balls is inversely proportional to the radius of the ball.
One says that there is an \emph{inverse  power-law relationship} between the two quantities. 
Similarly, if we consider a cube in $\mathbb{R}^3$, then we find that if we rescale the radius of balls needed to cover it by $1/s$, the minimum number of balls needed to cover it will change by $s^3$. The inverse power-law relationship between $N(\epsilon)$ and $1/\epsilon$ in this case is given by 

\[
\lim_{\epsilon\to 0}N(\epsilon)=\left (\frac{1}{\epsilon}\right )^3\, .
\]

\noindent The exponent in the inverse power-law relationship is equal to the box-counting dimension. In these two cases, these agree with the usual topological notion of dimension. 
The \emph{topological dimension} of a subset of Euclidean space can be defined recursively: it is $0$ if each point has arbitrarily small neighbourhoods each with empty boundary, and for $d\geq 1$ it is $d$ if each point has arbitrarily small neighbourhoods each with boundary of dimension $d-1$. In particular, the topological dimension of the empty set and of a finite set of points is $0$, while it is $1$ for a non-empty interval, and $3$ for a cube. More generally,   it is always a non-negative integer. %

\begin{figure}[t!]
\begin{center}
    \includegraphics[width = 0.35\textwidth]{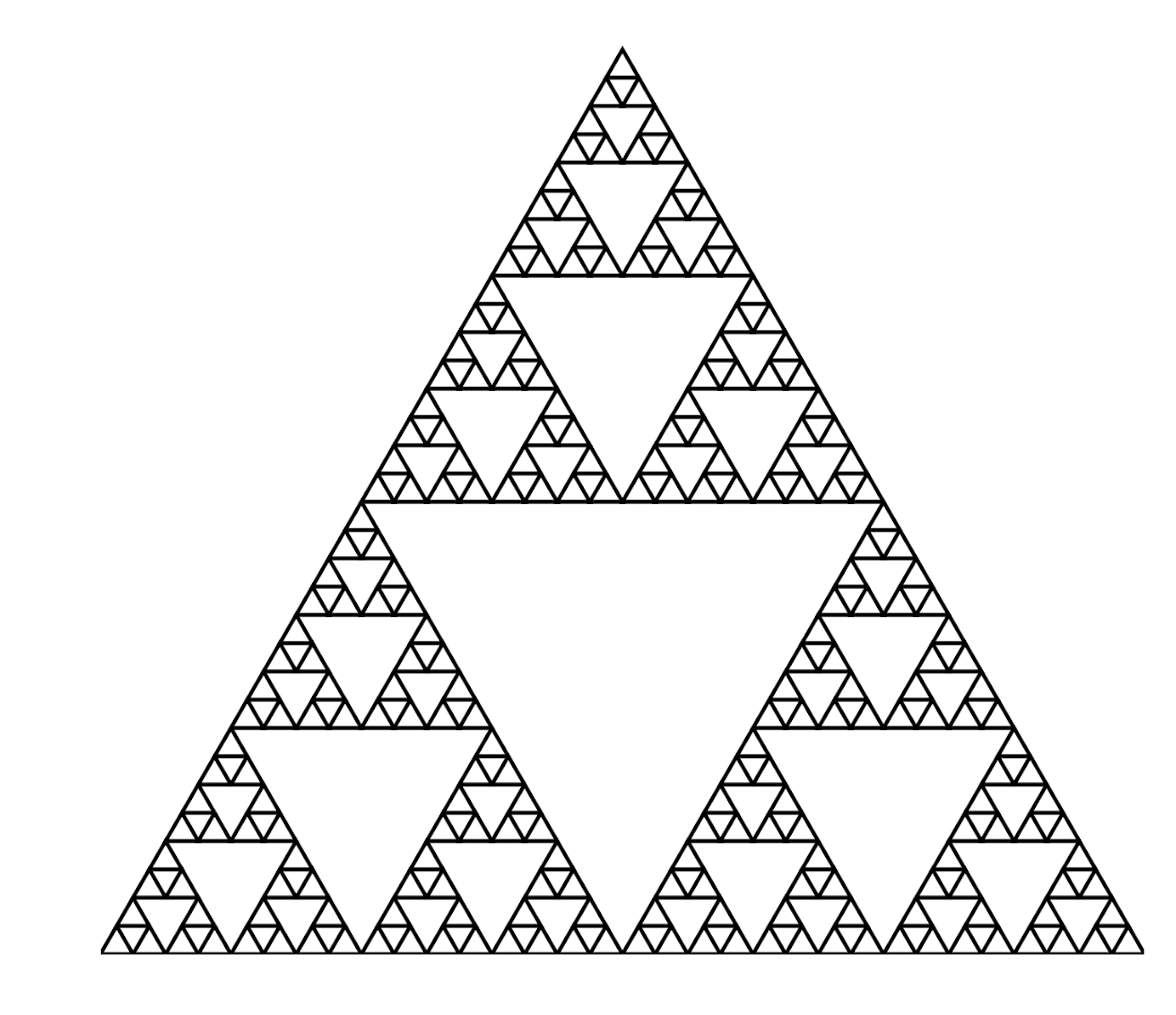}
    \end{center}
\caption{The  \emph{Sierspinki triangle},  a subset of $\mathbb{R}^2$ with box-counting dimension $\log(3)/\log(2)$.  Here we depict an approximation of the Sierpinski triangle, given by  the first 5 levels of an iterative construction of the set.}
\label{F: ex self-sim}
\end{figure}

One way to interpret the box-counting dimension  
is that it captures the self-similarity 
 of subsets of Euclidean space. Roughly,   a subset $S\subset \mathbb{R}^n$ is \emph{self-similar} if it can be decomposed into a number of disjoint pieces, each of which is a rescaled copy of $S$. An example of self-similar set studied in fractal geometry is the Sierpinski triangle $T\subset \mathbb{R}^2$, depicted  in Figure \ref{F: ex self-sim}: we can for instance decompose the Sierpisnki triangle into three disjoint subsets, each of which is a copy of the original set, rescaled by $1/2$.  A precise definition of self-similar sets can be given in terms of iterated function systems, see \cite[Chapter 9]{falconer2004fractal}. While we do not provide such a definition here, for the Sierpisnki triangle one can give a recursive definition as follows. At the $0$th level of iteration one considers an equilateral triangle in $\mathbb{R}^2$ with side-length $1$. At   level $i$, for $i\geq 1$, three copies of the previous space, rescaled by $1/2$ and arranged on the three corners of the original triangle. In Figure \ref{F: ex self-sim} we depict the  $5$th level of this iterative construction. 
To compute the box-counting dimension of the Sierpinski triangle, we note that if we rescale the radius of balls that cover it by $1/2$, we will need three times as many balls to cover it. More generally, in  this case we obtain the relationship 
\[
\dim_{box}(T)=\lim_{\epsilon\to 0}\frac{\log(N(\epsilon))}{\log(1/\epsilon)}=\frac{\log(3)}{\log(2)} \, .
\]

On the other hand, the topological dimension of $T$ can be shown to be equal to $1$.
More generally, one always has that for any $S\subset \mathbb{R}^n$  its box-counting dimension is greater or equal to its topological dimension.

In the classical monograph by Mandelbrot \cite{mandelbrot1983fractal}, a fractal was defined to be a subset of Euclidean space whose box-counting dimension is strictly greater than its topological dimension.\footnote{More precisely, Mandelbrot's definition was given in terms of the Hausdorff dimension. The notion of  \emph{Hausdorff dimension} is  more subtle;  roughly, it can be understood as studying the growth of the number of sets in covers of a given diameter --- in particular,  the sets in the cover do not all need to have the same diameter, as is the case in the box-counting dimension.
Box-counting and Hausdorff dimension are not equivalent in general, but they agree on a large class of spaces \cite{Bishop_Peres_2016}. } However, today this definition is understood as being too restrictive. While no agreed-upon definition of fractal exists, there is rather an agreement on properties shared by ``fractal-like'' subsets of Euclidean space, including that (i) $S$ is fractal if it has fractal dimension greater than its topological dimension, and that  (ii) if  $S$ is fractal then it contains  irregularities/non-smoothness at all scales.

In the current article, when referring to ``fractals'', we thus do not refer to a precise definition, but rather refer to spaces that satisfy such  properties. 
We note that in practice, not all fractals are self-similar. There might be fractals that are approximately so or not at all. And conversely, not all self-similar spaces are fractal, as illustrated for instance by a line in $\mathbb{R}^2$ which is self-similar. However, the intuition behind the definition  of  dimension of self-similar spaces can help in defining and computing the dimension of more general subsets of Euclidean space.
While here we choose to introduce the box-counting dimension, there are many different notions of dimension that generalise the topological notion of dimension to beyond the integers, including the Hausdorff dimension, which is one of the most well-known ones. These different notions are inequivalent in general, but they agree on certain classes of spaces. We will refer to such different notions collectively as ``fractal dimensions''.   We refer the reader to the monographs \cite{falconer2004fractal,mandelbrot1983fractal} for further details on the mathematical theory of  fractals.

\subsection{Fractal dimensions for networks: current approaches}
The interpretation of box-counting dimension of subsets of Euclidean space as being able to capture their self-similarity indicates a possible way to generalise such a dimension to weighted networks. A valid candidate for  dimension for  networks ought to measure the repeating of patterns in subnetworks of different sizes.
We will devote the rest of this section to making this statement precise, and to giving a summary of the current approaches.

\subsubsection{Self-similar networks}

Many physical or natural phenomena that are usually modelled with networks exhibit fractal-like properties  --- i.e., irregularities or  patterns repeating across subnetworks of different sizes. Some examples include tracheobronchial trees in humans \cite{West87}, the growth of blood vessels \cite{GBL95},  the world wide web \cite{albert2002statistical}, social network dynamics \cite{liu2016self}, 
epidemic propagation risk \cite{NIAN24}, 
or  the network
of neurons in the brain \cite{sporns}. 
Similarly as for self-similar subsets of Euclidean space, one may give a recursive definition of self-similar networks \cite{CARLETTI20102134}, as we explain next.

Let $G_0=(V,E,w)$ be a weighted network, i.e., we have a set $V$ of nodes, a set of edges $E=\{\{u,v\}\subset \mathcal{P}(V)\mid u\ne v\}$ together with a weighting function $w\colon E\to \mathbb{R}$.
Let $f\in \mathbb{R}$ and $s\in \mathbb{N}$ and assume that $0<f<1$ and $s>1$.
 Let $a\in V$ be a choice of node. 
 Consider $s$ copies $G_0^{(1)},...,G_0^{(s)}$ of $G_0$, whose weighted edges have all been scaled by the factor $f$.
 For $i = 1,...,s$ let us denote by $a^{(i)}$ the node in $G_0^{(i)}$ corresponding to a copy of the node $a \in G_0$. Now,  we define $T_{s,f,a}(G_0)$ to be the network obtained by joining each   $a^{(i)}$ to an additional node $a'$ through an edge of weight $1$.

Thus, we may define recursively 
$G_k=T_{s,f,a}(G_{k-1})$ for any $k\geq 1$, where we call each iteration $k$ a \emph{level}. In particular, the  network $G_\infty$ may be defined as the fix-point of the map $T_{s,f,a}(G_\infty)=G_\infty$, see \cite{CARLETTI20102134} for further details. 

\begin{figure}[h!]
\begin{center}
\includegraphics[width = 0.5\textwidth]{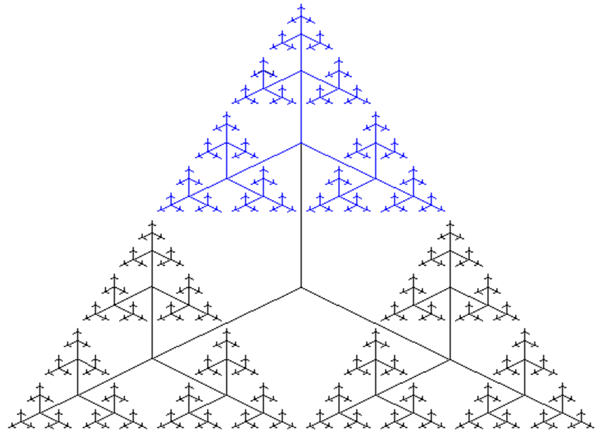}
\end{center}
\caption{Example of self-similar network. A network analogue of the   Sierpinski triangle, called \emph{Sierpinski tree}, which exhibits repeating patterns at all scales: if we consider a subnetwork of a given diameter $r$ (e.g., the one depicted in  blue) and we rescale its diameter by $2$, we obtain a subnetwork that contains  $3$ copies of the original subnetwork (e.g., for the blue subnetwork we obtain  the whole network depicted here). The Sierpinski tree is an infinite network, in this illustration we depict the first $5$ levels of the network of its recursive construction (i.e., $G_5$).}\label{F: ex sierp tree}
\end{figure}

We provide an example of a  network $G_\infty$ constructed through this recursive procedure in Figure \ref{F: ex sierp tree}: in this case $G_0$ is the network consisting of one node and no edges, and the parameter values are $s=3$ and $f=1/2$. Since we cannot depict an infinite network, in the illustration we represent the network  $G_5$ given by the fifth level of the recursive construction. 

One can thus give the following definition:
\begin{definition}\label{D:self-sim netw}
A weighted network $G=(V,E,w)$ is \define{self-similar} if there exists $a\in V$, $s\in \mathbb{N}_{\geq 1}$,
$f\in \mathbb{R}$ with $0<f<1$ such that $T_{s,f,a}(G)=G$.
\end{definition}
We note that, similarly as for subsets of Euclidean space, not every self-similar network is  ``fractal-like'', in the sense of containing ``irregularities'' at all scales (see also discussion at the end of Subsection \ref{SS:box-counting}). For instance, the line network, i.e., a network with  nodes  $V=\mathbb{N}$ and with edges $E=\{\{i,i+1\}\mid i\in \mathbb{N}\}$, all with weight $1$, is self-similar in the sense of Definition \ref{D:self-sim netw}. However, it does not contain irregularities at all scales, since any connected subnetwork is a copy of the whole network.  

The usefulness of considering self-similar networks is that any notion of dimension that generalises the box-counting dimension in an appropriate sense ought to give a dimension equal to  $\log(s)/\log(1/f)$ for such a network. Thus, they can serve as a testing ground for new notions of dimensions for networks.

\subsubsection{Fractal dimensions of weighted networks}
Many of the current  approaches to define fractal dimensions for weighted networks rely on finding  appropriate ways to cover the network, and studying the growth of the number of sets in a minimal cover as one rescales the sizes of the sets in the cover. The analogue of a covering ball or box in a network is a subnetwork, i.e., a subset of nodes together with all edges between the nodes that are present in the original network.
The analogue of the radius of a ball is given by the diameter of such a subnetwork, which in turn relies on a canonical metric structure on weighted networks.

\begin{definition}
Let $G=(V,E,w)$ be a weighted network. The \define{shortest-path metric} $\SP$  is a metric defined on the set of vertices $V$ as follows:
\begin{alignat}{2}
\SP\colon V\times V&\notag\to [0,\infty]\\
(u,v) &\notag\mapsto \min_{\gamma} \length(\gamma)
\end{alignat}
where the minimum is taken over all edge paths $\gamma=(e_1,\dots, e_m)$ connecting two nodes $u$ and $v$, and the length of such a path is given by the sum of the weights of its edges, i.e., $\length(\gamma)=\sum_{i=1}^m w(e_i)$.

We define the \define{diameter} of a weighted network $G$ to be the diameter with respect to its shortest-path metric.

\end{definition}

An example of a cover considered is the following.

\begin{definition}
Let $G=(V,E,w)$ be a weighted network, together with the shortest-path metric $\SP$.  
Let $\epsilon>0$. An \define{$\epsilon$-node covering} for $G$ is a  collection of  subnetworks $G_i$ each with  diameter bounded by $\epsilon$ (where we consider the induced metric on the subnetworks), and such that  each node of the network $G$ is contained in exactly one subnetwork. 
\end{definition}

One can then define a network analogue of the box-counting dimension as follows:

\begin{definition}\label{D:box-counting netw}
Let  $G=(V,E,w)$ be a weighted network.
Let $N(\epsilon)$ be the number of subnetworks in a minimal $\epsilon$-node covering of the network. Then one defines
\[
\dim_{box}(G)\lim_{\epsilon\to 0}\frac{\log (N(\epsilon))}{\log(1/\epsilon)}
\]
to be the \define{box-counting dimension} of $G$. 
\end{definition}

There are many other types of coverings that have been studied, such as covering in which nodes can belong to more than one subnetwork \cite{SZ14}, or \emph{edge coverings}, in which one  requires each edge, instead of node, to belong to exactly one subnetwork \cite{ZJS07}, to give a few examples.

Other approaches to defining fractal dimensions for networks rely on studying neighbourhoods of nodes, instead of covering the network. 
\begin{definition} Let $G=(V,E,w)$ be a weighted network, and let $\SP$ be its shortest-path metric. 
     Given $\epsilon>0$ the \define{$\epsilon$-neighbourhood} of a node $x\in V$ is defined as follows:
     
     \[
     N(x,\epsilon)=\left \{ v\in V \mid \SP(v,x)\leq \epsilon\right \} \, .
     \]
     In other words, the neighborhood consists of  all the nodes at shortest-path distance  at most $\epsilon $ from $x$.
\end{definition}

\begin{definition}
Let $G=(V,E,w)$ be a weighted network. Let $x\in V$ be a node. Denote by $N(x,\epsilon)$ the $\epsilon$-neightbourhood of $x$. The \define{internal scaling dimension  starting from $x$} of $G$ is given by 
\[
\dim_{\mathrm{int}}(x)=\lim_{\epsilon\to \infty} \frac{|N(x,\epsilon)|}{\log(\epsilon)}
\]
if the limit exists, where we use $|S|$ to denote the cardinality of a given set $S$.

If for every $x\in V$ $\dim_{\mathrm{int}}(x)$ exists and if 
$\dim_{\mathrm{int}}(x)=\dim_{\mathrm{int}}(x') =:D$ for all $x\ne x'$ then one says that $G$ has \define{internal scaling dimension} $\dim_{\mathrm{int}}(x)=D$.
\end{definition}

Other methods relying on neighbourhoods to define fractal dimensions include the correlation dimension \cite{lacasa2013correlation}, and  surface dimension \cite{NR1998}. 
One particular advantage of  approaches based on neighbourhoods is that they can also be used to estimate fractal dimensions of networks with infinite diameter. Another approach to estimate fractal dimensions of networks with infinite diameter is mass dimension \cite{zhang2008transition}, which relies on estimating the diameter of a $1$-parameter family of subnetworks converging to the full network.

 For an overview of fractal dimensions for networks, we refer the reader to  the monograph \cite{rosenberg} or the survey article \cite{WEN202187}.

As we have discussed in this section, there is a plethora of  different approaches that have attempted to calculate various types of fractal dimension for networks. This multitude of approaches and resulting lack of consensus on notions of dimensions are due to the fact that transferring continuous methods like the box-counting method to discrete objects such as networks is challenging and as a consequence there is no consensus on a natural way to calculate these dimensions for networks. We believe that topological methods might provide an effective method for estimating this dimension, 
as we explain next.

\section{Existing topological approaches to estimating fractal dimensions}

Topological Data Analysis (TDA) \cite{carlsson09} is a field that uses methods and insights from  Topology --- the area of mathematics that studies ``abstract shapes''  --- to study local as well as global properties of data, with a variety of applications, including hierarchical clustering \cite{clustering2007,tomato}, dimensionality reduction \cite{bei2018,dreimac2023}, signal denoising \cite{imdenoising} and analysis \cite{perea}, complex networks \cite{petri2013topological,horak2009persistent}, protein folding \cite{prfolding,xia2014persistent}, the analysis of graphs \cite{carriere2020perslay}, 3D objects \cite{3Dchazal,Reininghaus_2015_CVPR,turner2014persistent}, 2D images \cite{carlsson2008local,mnist}, or heterogeneous and high-dimensional datasets \cite{hetdataset,carlssondatasets}. For a  database of TDA applications to real-world data sets, see \cite{DONUT}.

\subsection{Persistent homology}

One of the most popular methods in TDA is, arguably,   persistent homology (PH) \cite{ZC05}.  %
In a nutshell, persistent homology-based methods associate to a data set, such as a finite metric space or a network,  
 a simplicial complex along with a filtration (a way to build it as a sequence of nested sub-complexes).
 Then, one can detect and count the number of $i$-dimensional holes --- components ($i=0$), holes ($i=1$), voids ($i=2$ and higher-dimensional voids ($i\geq 3$)) --- , track their emergence and disappearance  as the filtration parameter varies, and store this information in compact summaries called \emph{persistence barcodes}. Notably, both local and global properties, such as curvature and convexity \cite{bubenik2020persistent,turkes2022effectiveness} may be extracted. One of the most useful properties of persistent homology is stability: small variations in input data result in persistence barcodes which are close each other \cite{cohen2005stability, wasserstein-stability}. Another key feature of PH is that there are many choices for associating a simplicial complex and filtration to a data set, hence, allowing one to work with different types of data, such as networks.

 Here we provide  some key definitions and examples illustrating persistent homology computations on finite metric spaces, which  will be needed to illustrate PH-based notions of fractal dimensions for subsets of Euclidean space.

\begin{definition}[Vietoris-Rips Complex]\label{def:VR_complex}
    Let $(X,d) \subset \mathbb{R}^D$ be a finite metric space with $n \in \mathbb{N}$ points. 
    Given a parameter $\epsilon > 0$, the \define{Vietoris-Rips (VR) complex of $X$ at scale $\epsilon$} is the (abstract) simplicial complex defined as follows:
    
    \[
    \VR_\epsilon(X)=\{\{x_1,\dots , x_n\}\in \mathcal{P}(X) \mid d(x_i,x_j)\leq \epsilon \forall i,j\}\, ,
    \]
    where $\mathcal{P}(X)$ denotes the power set of $X$.
\end{definition}

Given $\epsilon\leq \epsilon'$ we note that $\VR_{\epsilon}(X)\subseteq \VR_{\epsilon'}(X)$. More generally, we give the following definition:

\begin{definition}
Let $K$ be an abstract simplicial complex. A \define{filtration} of $K$ is a collection 
     of (abstract) simplicial complexes 
$\left \{K_\epsilon\right\}_{\epsilon\in \mathbb{R}}$ 
such that $K_{\epsilon}(X)\subseteq K_{\epsilon'}(X)$ whenever $\epsilon\leq \epsilon'$ and $\cup_\epsilon K_\epsilon = K$. 

We call a simplicial complex $K$ together with a given filtration a \define{filtered simplicial complex}.
\end{definition}

Computing simplicial homology in a given degree $i$ with coefficients in a field $\mathbb{K}$ of a filtered simplicial complex yields what are called ``persistence modules'' in the TDA literature:

\begin{definition}\label{D:pers mod}
Let $\mathbb{K}$ be a field. A \define{persistence module} over $\mathbb{K}$ is given by a tuple 
$M={\left (\{V_\epsilon\}_{\epsilon\in \mathbb{R}},\{\phi_{\epsilon,\epsilon'}\}_{\epsilon\leq \epsilon'}
\right )}$
where each $V_i$ is a $\mathbb{K}$-vector space, and each  $\phi_{\epsilon,\epsilon'}\colon V_\epsilon\to V_\epsilon'$ is a $\mathbb{K}$-linear map satisfying the properties that 
$\phi_{\epsilon',\epsilon''} \circ\phi_{\epsilon,\epsilon'}$ whenever $\epsilon\leq \epsilon'\leq \epsilon''$  and further  $\phi_{\epsilon.\epsilon}=id_{V_\epsilon}$ for each $\epsilon\in \mathbb{R}$.
\end{definition}

For instance, given a filtration of VR complexes $\{\VR_\epsilon(X)\}_{\epsilon\geq 0}$, we obtain one persistence module $\left (\{H_i(\VR_\epsilon(X))\}_{\epsilon},\{\phi_{\epsilon,\epsilon'}\}_{\epsilon\leq \epsilon'} \right )$ for each homology degree $i$, where the linear maps $\phi_{\epsilon,\epsilon'}$ are the maps induced by the inclusions $\VR_{\epsilon}(X)\subseteq \VR_{\epsilon'}(X)$.

\begin{figure}[ht!]
    \centering
    \includegraphics[width=1\textwidth]{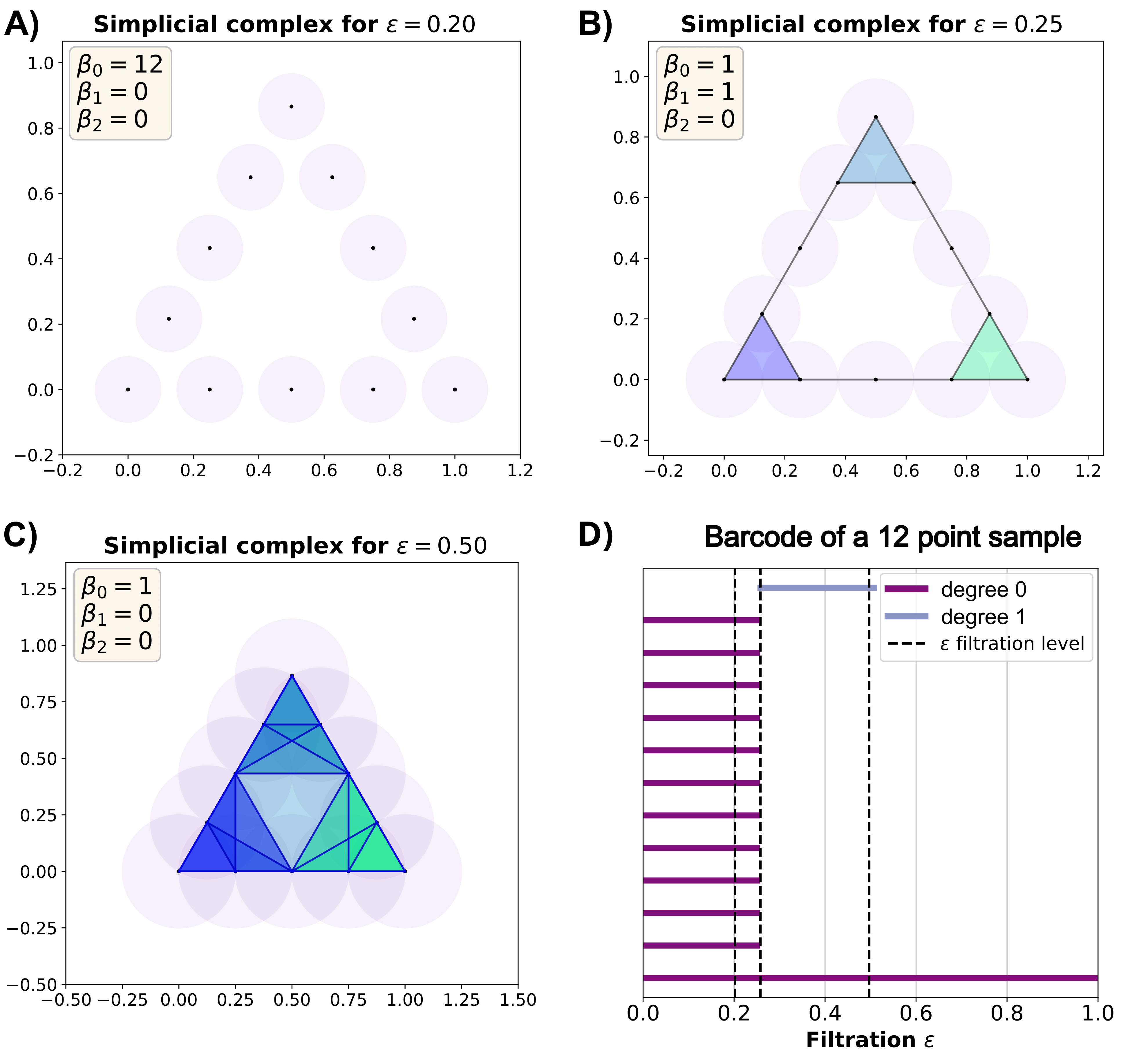}
   \caption{A)-C): Examples of VR complex on a subsample of 12 points from a triangle in $\mathbb{R}^2$,
    shown for different filtration parameters. We indicate the ranks $
    \beta_i$ of the homology vector spaces in degree $i$ for each filtration value.  D): The resulting barcode of the full filtration, with dotted lines highlighting the chosen filtration values in A)-C) for homology degrees $0$ and $1$.}
    \label{ex:persistence_homology_triangle}
    \end{figure}

It is a fundamental theorem in TDA that one can summarise not only the topological properties --- i.e., the number of $i$-dimensional holes, given by the ranks $\beta_i$ of homology vectors spaces in degree $i$  ---   of each simplicial complex in a given filtration, but also  keep track of how these properties change as one varies the filtration parameter --- whether two components merge, holes or voids appear or are filled in --- through an algebraic invariant of the persistence module of the filtration, called  \emph{persistence barcode}, one for each homology degree $i$, denoted by $B_i(X)$. Roughly, such a summary consists of a collection of disjoint intervals, each with left-endpoint corresponding to the appearance, or \emph{birth}, of a given $i$-dimensional hole and right-endpoint corresponding to the merging or disappearence, the \emph{death}, of the feature. We refer the reader to \cite{otter2017roadmap} for a friendly introduction to the theory and computation of persistent homology.
 We provide an example of persistence barcode in  Figure ~\ref{ex:persistence_homology_triangle}.  %

\begin{example}[VR filtration on points sampled from a  triangle]\label{ex:PH_on_Sierpinski_triangle}
In Figure \ref{ex:persistence_homology_triangle}, we show an example of VR filtration on a finite metric space consisting of 12 points sampled from a triangle in $\mathbb{R}^2$. 
   
\end{example}

We conclude this subsection with two further examples of simplicial complexes widely used in TDA, which will be relevant for the definition of persistent homology dimension, and persistent magnitude dimension, respectively.
\begin{definition}[\v{C}ech Complex]\label{def:cech_complex}
    Let $(X,d) \subset \mathbb{R}^D$ be a finite metric space with $n \in \mathbb{N}$ points. 
    Given  $\epsilon > 0$, consider the cover of $X$ consisting of closed balls $\{B_\epsilon(x)\}_{x\in X}$ of radius $\epsilon$ centered at the points of $X$ . The \define{\v{C}ech complex of $X$ at scale  $\epsilon$} is the  (abstract) simplicial complex 
    defined as follows:

    \[
\mathcal{C}_\epsilon(X)=\left\{\{x_1,\dots , x_n\}\in \mathcal{P}(X) \mid \bigcap_{i=1}^nB_\epsilon (x_i)\ne \emptyset
\}\right\}\, .
    \]
\end{definition}
While the \v{C}ech complex and Vietoris-Rips complex constructions yields different simplicial complexes with different topological properties, there is a precise way in which one can see them as being approximations of each other, see \cite{edelsbrunner2010computational}.

\begin{definition}\label{D:alpha}
Let $(X,d) \subset \mathbb{R}^D$ be a finite metric space with $n \in \mathbb{N}$ points. 
For each $x\in X$, let $V_x$ denote the set of points in $\mathbb{R}^D$ that are closer to $x$ than any other point in $X$, namely
\[
V_x=\{y\in \mathbb{R}^D \mid d(x,y)\leq d(x',y) \forall x'\in X\}\, .
\]
The \define{alpha complex of X at scale $\epsilon$} is defined as the following abstract simplicial complex:

 \[
\mathcal{\alpha}_\epsilon(X)=\left\{\{x_1,\dots , x_n\}\in \mathcal{P}(X) \mid \bigcap_{i=1}^n\left (B_\epsilon (x_i)\cap V_{x_i}\right)\ne \emptyset
\}\right\}\, .
    \]
\end{definition}

While the alpha complex captures exactly the same topological information as the \v{C}ech complex (i.e., they have the same homotopy types, see \cite{edelsbrunner2010computational}), the alpha complex has the advantage that it can yield sparse simplicial complexes with fewer simplices in high dimension with respect to the \v{C}ech complex.

\subsection{Persistent homology and fractal dimensions}\label{SS:PH dim}

One of the earliest applications of topology to the study of data appeared in \cite{robins00}, 
in which it was found that the growth rates of the number of components or holes (i.e., \emph{persistent Betti numbers}) of certain subsets of Euclidean space are closely related to their Minkowski dimension. %
A number of different fractal dimensions based on PH have since been proposed in the literature. One of the earliest ones \cite{macpherson2012measuring} introduces a notion based on the distribution of intervals in a persistence barcode and demonstrates that it can be used to characterise self-similar structure in a number of examples  --- such as branched polymers, Brownian trees, and self-avoiding random walks ---, and that it coincides with the Hausdorff dimension, albeit for a limited number of cases. It has been %
proposed later in \cite{Jaquette2019FractalDE}, that this quantity, rather than a dimension, may be interpreted as an indicator of the difficulty in estimating  so-called PH-dimensions (see below) and the correlation dimension of a shape.%

Subsequently, two different notions of fractal dimensions, both  called \emph{PH-dimensions},  have been suggested in \cite{schweinhart2020fractal,adams2020fractal} for measures on bounded metric spaces. These are collections of dimensions $PH_i$, one for each homology degree $i\in \mathbb{N}$ (see Definition \ref{def:PH_dim}).
It has to be noted that, when both exist, one can recover the PH-dimension defined in \cite{adams2020fractal} [Definition~\ref{def:PH_dim}] as a special case of the  one  in \cite{schweinhart2020fractal}.
In \cite{schweinhart2020fractal}, the authors demonstrate that for a certain technical conditions on the measure (Ahlfors-regularity), their definition agrees with the Hausdorff dimension, which is the first result linking theoretically both dimensions. 
Here we introduce the more general notion from \cite{schweinhart2020fractal}.

Given a subsample of a  metric space $X$, the persistent homology dimension of  $X$ is defined using the decay rate of the length of short intervals in the persistence barcode of the subsample,  as the subsample size increases.
This decay rate is tracked using so-called \emph{power-weighted sums} (Definition \ref{def:power_weighted_sum}), namely weighted sums of lengths of intervals in a barcode, where different choices for the weighting parameter result in assigning different weight to intervals of different length (larger values of the weight parameter assign relatively more importance to longer intervals than shorter ones).

\begin{definition}[Power-weighted sum ]\label{def:power_weighted_sum}
Let $(X,d)$ be a metric space and let 
$S\subset X$ be a finite 
subset of points. 
  Let $\alpha >0$. %
 Let $B_i(S)$ be the persistence barcode in degree $i$ for
 a given filtered simplicial complex associated to $S$. The \define{$\alpha$-power-weighted sum} of $B_i(S)$ is given by 

    \[ E_{\alpha}^i(S) = \sum_{\underset{|I|<\infty}{I \in B_i(S)} }|I|^\alpha \, . \]
\end{definition}

We note that  the definition of power-weighted sum relies on a choice of filtered simplicial complex. In \cite{Jaquette2019FractalDE} the authors use  
 the \v{C}ech complex (see Definition\ref{def:cech_complex}), while in \cite{adams2020fractal} the Vietoris-Rips complex (see Definition \ref{def:VR_complex}) is chosen. 
 
 For the next definition we consider a fixed choice of  filtered simplicial complex over all subsamples. 
One then defines:

\begin{definition}[$i$-dimensional Persistent Homology Dimension]\label{def:PH_dim}
Let $(X,d)$ be a metric space, and let $\mu$ be a probability measure on it. 
Let $S_n\subset X$ be 
a  subsample of cardinality $n$ drawn according to the probability measure
$\mu$.
Let $\alpha > 0$.
Let
    \[ \beta = \lim \sup_{n \rightarrow \infty} \frac{\log(\mathbb{E}(E_{\alpha}^i(S_n)))}{\log(n)} \, . \]
   
Then the \define{$i$-dimensional PH dimension of $X$} is defined as

\[ \mathrm{dim}_{PH_i^\alpha}(X,\mu) = \frac{\alpha}{1-\beta}\,  . \]
\end{definition}

\begin{figure}
    \centering
    \includegraphics[scale=0.25]{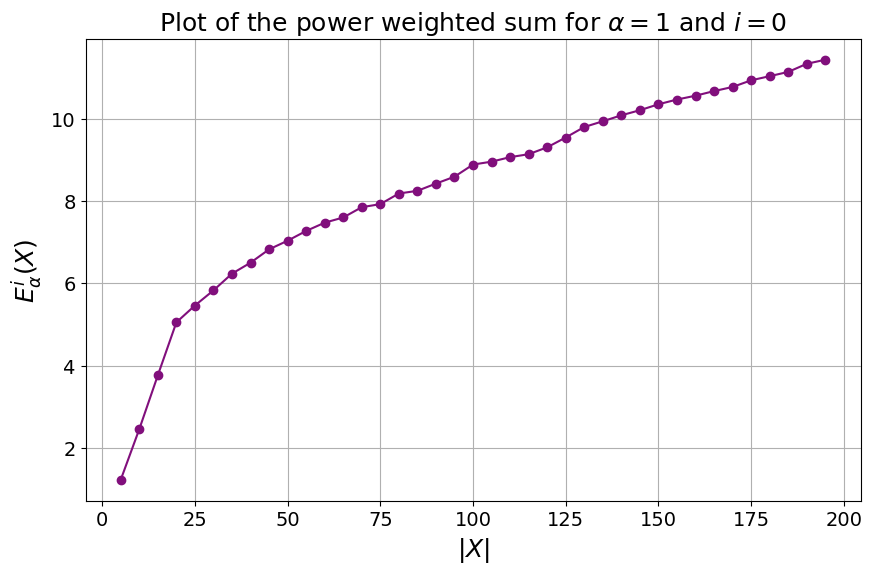}
    \includegraphics[scale=0.25]{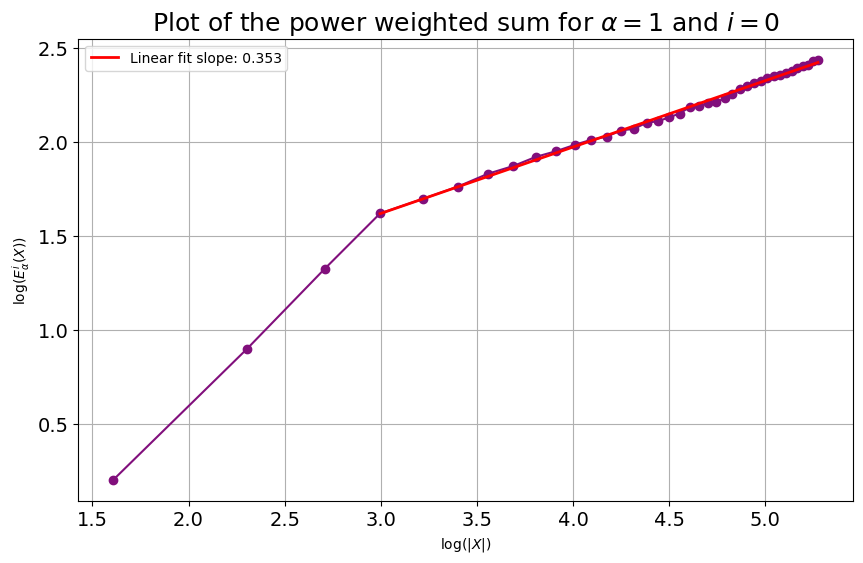}
    \caption{(Left)  Power-weighted sums for   different subsamples of the Sierpinski triangle at level $7$  plotted against the subsample sizes and (Right) the resulting log-log plot, together with a red line illustrating a linear regression fit.}
    \label{fig:log-log-plots}
\end{figure}

\begin{example}[PH-dimension of Sierpinski triangle] 
For this example, we 
consider an approximation of the Sierpinski triangle  at level $7$. Furthermore, we fix $\alpha = 1$. 
Consider subsamples of $n$ points in the Sierpinski triangle for $n = 5, 10,\dots, 200$ ($40$ subsamples total). In Figure \ref{fig:ex_Sierpinski_level7}, we show the VR filtration for different subsamples ($n=25,50,200$) of the Sierpinski triangle and the resulting persistence barcodes. One may observe that the shorter intervals in the barcodes in degree $0$ become shorter as one increases the subsample size. The rate at which their length decreases encapsulates information about the fractal dimension of the shape. In Figure \ref{fig:log-log-plots} we show the resulting graph of the power-weighted sums against the subsample sizes and the resulting log-log plot, respectively. To compute the estimated growth rate $\beta$, we fit a linear regression to the last $36$ points in the log-log plot. We obtain $\beta \approx 0.353$.

Then the 0-dimensional PH-dimension closely approximates the box-counting dimension ($1.58$) of the Sierpinski Triangle \cite{Jaquette2019FractalDE}:

\[\dim_{PH_0^1}(S) \approx 1.546. \]

\begin{figure}
    \centering
    \includegraphics[scale=0.27]{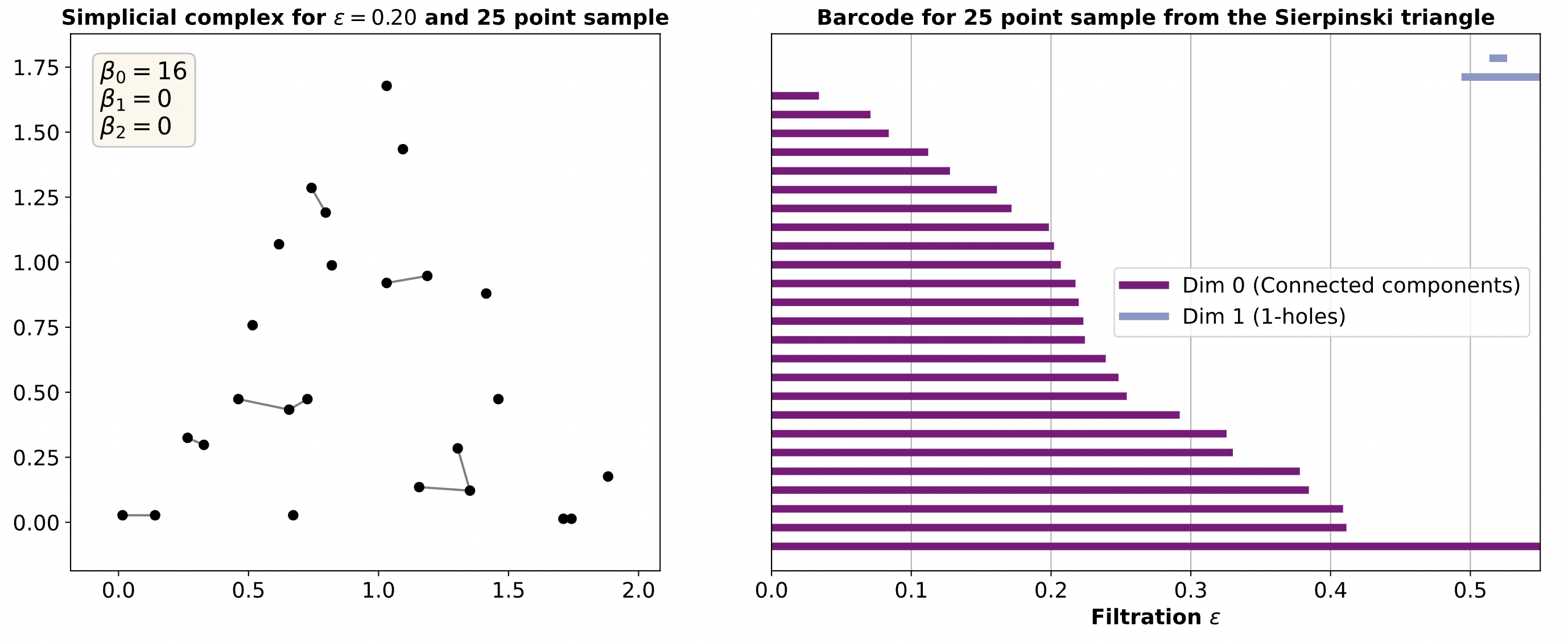}
    \includegraphics[scale=0.27]{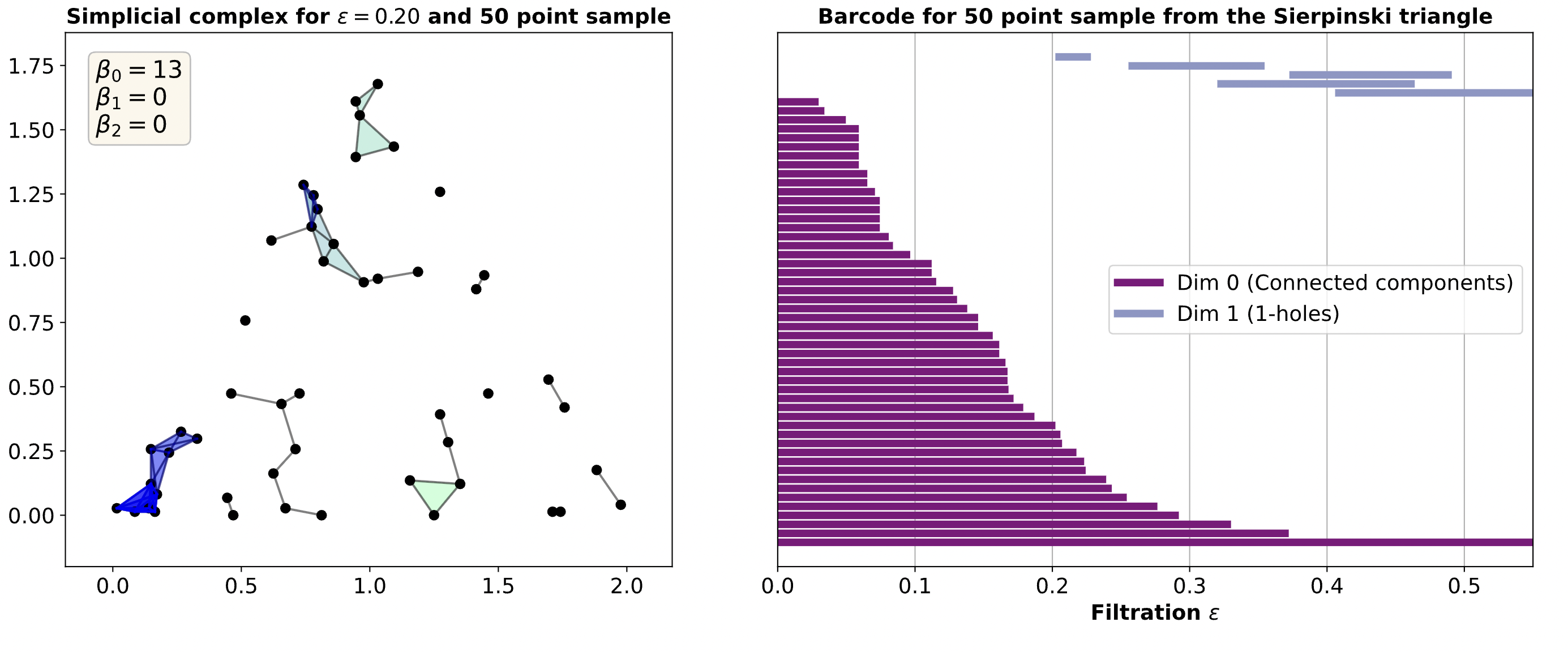}
    \includegraphics[scale=0.27]{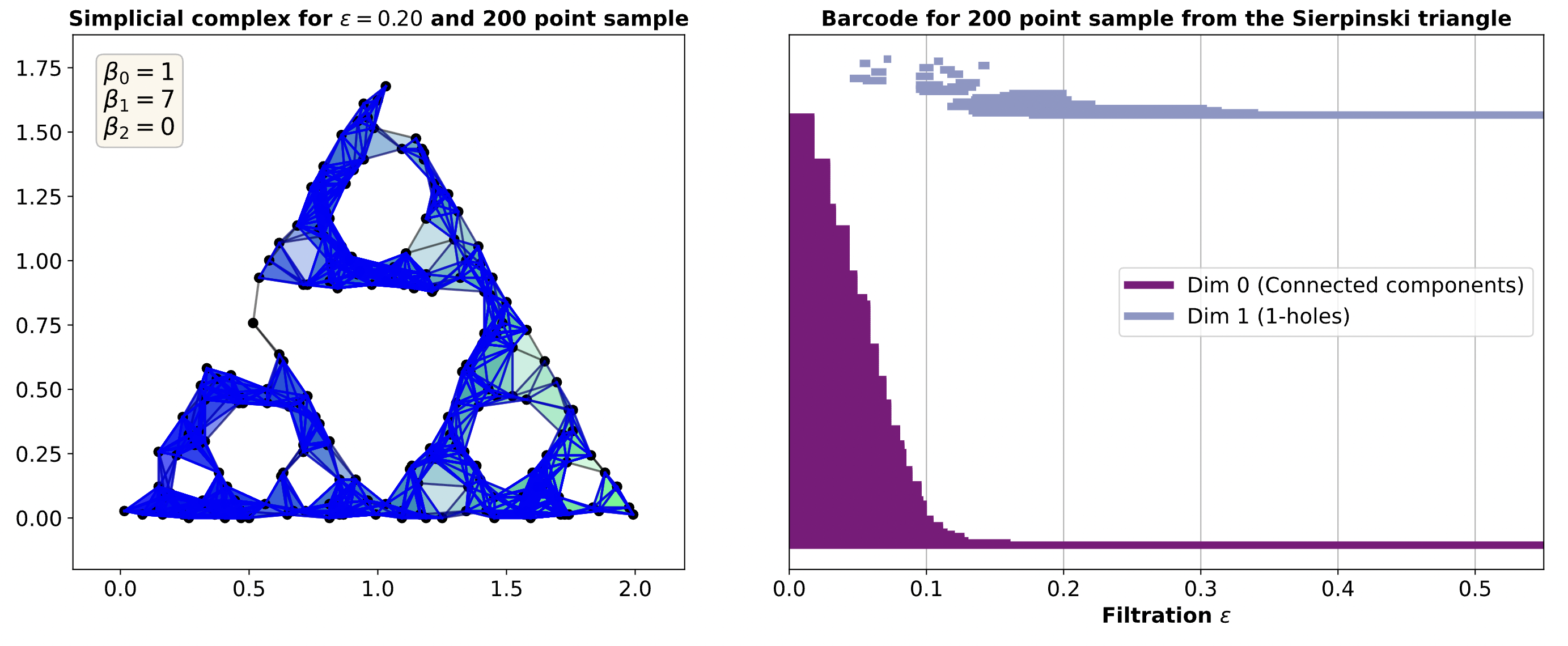}
    \caption{(Left) Vietoris-Rips complexes $\VR_{0.2}(X)$  for subsamples $X$ of $25$ (top), $50$ (middle) and $200$ (bottom) points from the Sierpinski triangle at level $7$ and (Right) barcodes computed from the filtered Vietoris-Rips complexes associated to the subsamples of points.}
    \label{fig:ex_Sierpinski_level7}
\end{figure}
\end{example}

One can then ask the question: how well do these PH-based notions of fractal dimensions compare with other dimensions? In \cite{Jaquette2019FractalDE}, an extensive practical comparison between classical and PH-based dimensions has been carried out %
for different types of well-known subset of Euclidean space that exhibit fractal-like properties (``classical fractals''), as well as for  strange attractors. In the case of classical fractals, the $PH_{0}$ dimension and the correlation dimension yield similar performance and converge to the true  dimension, whereas, in most cases, the box-counting, $PH_{1}$, and $PH_{2}$ dimensions perform worse. 
For the case of strange attractors, the results often widely differ for the different notions of dimension considered. However, since the true dimension of these spaces is not known, it is difficult to assess which notion performs best.

\subsection{Magnitude}\label{SS: magnitude}

In a parallel line of work, another family of fractal dimensions has been developed, similar in flavour to the PH-based dimensions, and based on a relatively new but rich geometric invariant of metric spaces called \emph{magnitude} \cite{leinster2013magnitude} with known links to TDA \cite{otter2021magnitude}.

\begin{definition}
    Let $(X,d)$ be a finite metric space. Then, the \define{similarity matrix} of $X$ is defined as the matrix with entries $\zeta_{xy} = e^{-d(x,y)}$ for $x,y\in X$.
\end{definition}

\begin{definition}[Magnitude]
Let $(X,d)$ be a finite metric space with similarity matrix $\zeta_{ij}$. 
If $\zeta_{ij}$ is invertible, the \define{magnitude} of $X$ is defined as
\begin{equation}
\mathrm{Mag}(X) = \sum_{ij}(\zeta^{-1})_{ij}.
\end{equation}

\label{magnitude}
\end{definition}

When $X$ is a finite subset of $\mathbb{R}^n$, then $\zeta_{ij}$ is a symmetric positive definite matrix as proven in \cite{leinster2013magnitude} (Theorem 2.5.3). Then, $(\zeta^{-1})_{ij}$ exists and, hence, magnitude exists as well.

Magnitude is best illustrated when considering a few examples of spaces with a small number of points. We first compute the magnitude of a $1$-point space, followed by the $2$-point space.

\begin{example}
    Let $X$ denote the metric space with a single point $a$. Then, $\zeta_X$ is a $1 \times 1$ matrix with $\zeta_X^{-1} = 1$ and using the formula for magnitude, we get $\mathrm{Mag}(X) = 1$.
\end{example}

\begin{example}\label{ex:2points}
Consider the space of two points. Let $X=\{a,b\}$ be a finite metric space where $d(a,b)=t$. Then
\[
   \zeta_X=\begin{bmatrix} 1 & e^{-t}\\ e^{-t} & 1\end{bmatrix},  
\]

\noindent so that $
    \zeta_X^{-1} = \frac{1}{1 - e^{-2t}} \begin{bmatrix} 1 & -e^{-t}\\ -e^{-t} & 1\end{bmatrix},$ and therefore 
\[
    \mathrm{Mag}(X)=\dfrac{2 - 2e^{-t}}{1 - e^{-2t}}=\dfrac{2}{1+e^{-t}}. \label{eq:1}
\]
\end{example}

\begin{definition}
Let $(X,d)$ be a finite metric space. We define $(tX, d_t)$ to be the metric space with the same points as $X$ and the metric $d_t(x,y) = td(x,y)$.
\end{definition}

\begin{definition}[Magnitude function]
The \define{magnitude function} of a finite metric space $(X,d)$ is the partially defined function 
\begin{alignat}{2}
(0,\infty)&\notag\to \mathbb{R}\\
t &\notag\mapsto \mathrm{Mag}(tX)\, .
\end{alignat}
\end{definition}

In particular, one can observe that the magnitude function for the $2$-point space from Example \ref{ex:2points} has value $1$ at $t=0$ and asymptotically approaches $2$ for $t\to \infty$. One can interpret this as telling us that the closer the two points are to each, the more the metric space looks like a $1$-point space, while the further they are apart, the more it looks like a $2$-point space. More generally, there is a precise sense in which one can interpret the magnitude function as capturing the ``effective number of points'' of a space, see, e.g.,  \cite{LeinsterMeckes}.

\subsection{Magnitude and fractal dimensions}

A notion of dimension, the \emph{magnitude dimension}, was given using the magnitude function, i.e., investigating its asymptotic growth; empirical evidence that this dimension approximates the Hausdorff dimension was provided in \cite{willerton2009heuristic}.

The notion of magnitude can be extended to compact subsets of Euclidean space,  by taking the supremum over the magnitudes of finite subspaces, see \cite{meckes2015magnitude} for a precise definition.  
One then can define magnitude dimension as follows:

\begin{definition}[Magnitude dimension]
When
\begin{equation}
\label{mag_dimension}
\mathrm{dim_{Mag}}X = \lim_{t \to \infty}\frac{\log(\mathrm{Mag}(tX)))}{\log t}   
\end{equation}
exists, we define this number to be the \define{magnitude dimension} of $X$~\cite{meckes2015magnitude}.
\label{magnitude_dimension}
\end{definition}
Also, in \cite{meckes2015magnitude}, it is shown that for $X\subset \mathbb{R}^D$   magnitude dimension turns out to be the same as the box-counting dimension.

\begin{example}[Magnitude dimension of the Sierpinski triangle]
We consider $X$  a set of points sampled from the Sierpinski triangle at level $7$. In Figure \ref{fig:magnitude-dimension-example} we plot the magnitude dimension against different rescaling values $t={1, 2, ..., 300}$, and by fitting a line on the resulting log-log plot we obtain the following estimation:

\[ \mathrm{dim_{Mag}}X = \lim_{t \to \infty}\frac{\log(\mathrm{Mag}(tX)))}{\log t} \approx \frac{y_{80} - y_{40}}{l_{80}-l_{40}}= 1.55.\]

\begin{figure}
    \centering
    \includegraphics[scale = 0.25]{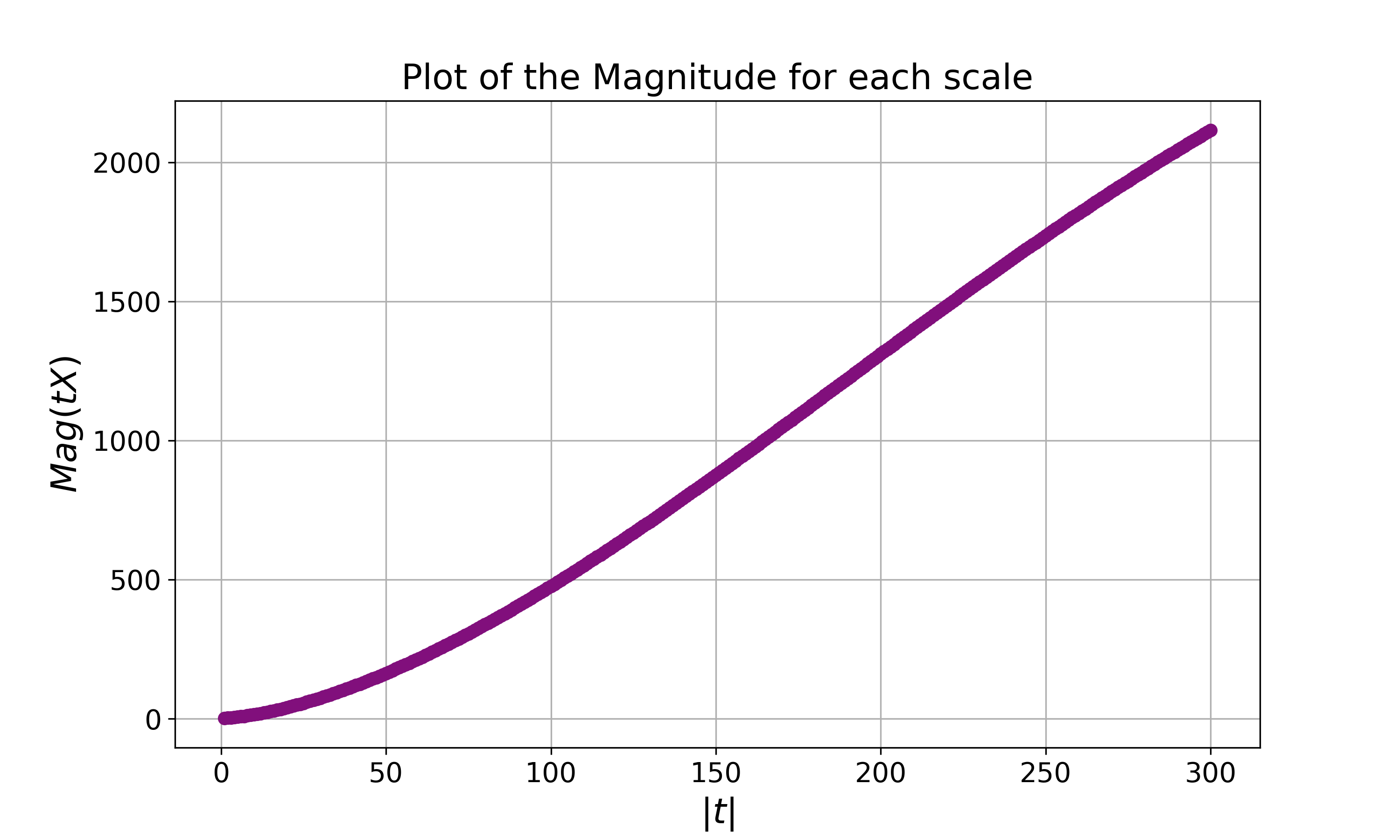}
    \includegraphics[scale = 0.25]{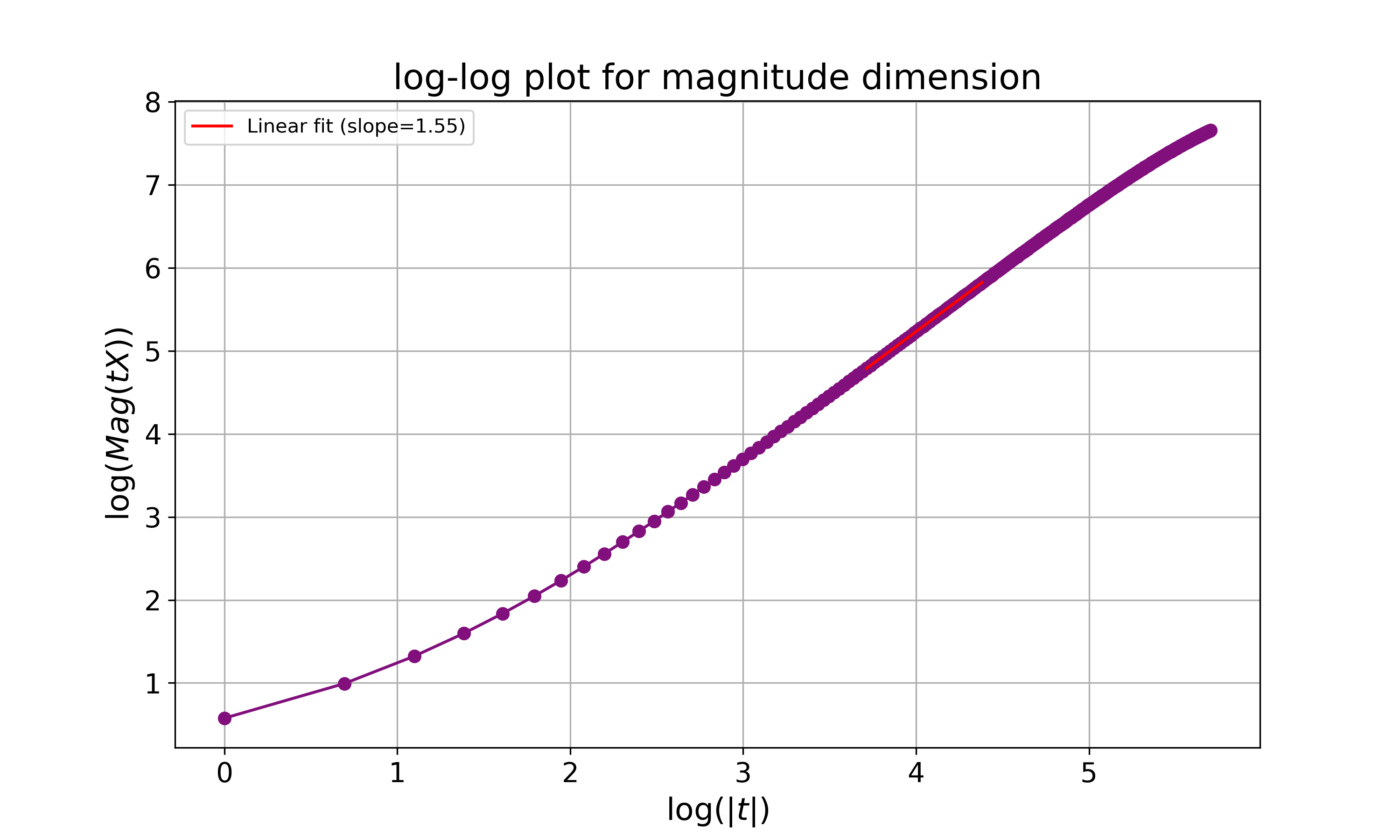}
    \caption{The  graph of the magnitude of the Sierpinski triangle against increasing scales (top). The resulting log-log plot used to estimate the magnitude dimension (bottom).}
    \label{fig:magnitude-dimension-example}
\end{figure}
   
 \end{example}

\subsection{Persistent magnitude dimensions}
Exploiting the link provided between the two fields in \cite{otter2021magnitude}, in \cite{govc2021persistent}, mag-
nitude and persistent homology are merged to introduce \emph{persistent magnitude},
a numerical invariant of graded persistence modules, and Rips magnitude as
a particular instantiation. 
It was formally introduced by Govc and Hepworth in \cite{govc2021persistent}, as a weighted and signed count of the intervals of the persistence module, in which an interval  of the form $[a, b]$ in homology degree $i$ is counted with weight $(e^{-a} -e^{-b})$ and sign $(-1)^i$.

For this, they consider \emph{graded persistence modules}: while computing simplicial homology in degree $i$ of a given filtered simplicial complex yields a persistence module $M_i$ (see Definition \ref{D:pers mod}), one could also consider all the information given by all homology degrees at once, and instead consider the $\mathbb{N}$-graded persistence module $M:=\oplus_{i\in \mathbb{N}} M_i$.

\begin{definition}[Persistent magnitude]
    Let $M$ be a finitely presented graded persistence module with barcode decomposition in degree $i$ given by  $\{[a_{i,j},b_{i,j}]\}_{j=1}^{m_i}$. The \define{persistent magnitude} of $M$ is the real number
    $$\mathrm{Mag}(M) = \sum_{i=0}^\infty \sum_{j=0}^{m_i} (-1)^i (e^{-a_{i,j}}-e^{-b_{i,j}}).$$
\end{definition}
This definition was inspired by the notion of \emph{blurred magnitude homology} \cite{otter2021magnitude}, which is the homology of a filtered simplicial set called the \emph{enriched nerve}.
In \cite{govc2021persistent} (Theorem 4.9) the authors showed that for any finite metric space $(X,d)$, the magnitude of the rescaled space $tX$ can be recovered from the barcode of the blurred magnitude homology of $tX$. More precisely, they show that for $t$ large enough one has
\[
\mathrm{Mag}(tX) = \sum_{i=0}^\infty \sum_{j=0}^{m_i} (-1)^i \left(e^{-a_{i,j}t}-e^{-b_{i,j}t}\right)\, ,
\]
where $\{(a_{i,j},b_{i,j})\}_{j=1}^{m_i}$ is the barcode of the homology in degree $i$ of the enriched nerve of $X$.

Govc and Hepworth showed that persistent magnitude has good formal properties, such as additivity with respect to exact sequences, and investigated the persistent magnitude for specific case of the Vietoris-Rips complex associated to a finite metric space, which they called \emph{Rips magnitude} of $X$, denoted by $\mathrm{Mag}_\VR(X)$.
In this case, one can observe that 
for a metric space $X$ of cardinality $n$  the infinite sum reduces to the following finite sum:
\begin{equation}\label{E:rips mag}
\mathrm{Mag}_\VR(tX) = \sum_{i=0}^{n-2} \sum_{j=0}^{m_i} (-1)^i \left(e^{-a_{i,j}t}-e^{-b_{i,j}t}\right)\, ,
\end{equation}
since the maximum dimension of a Vietoris-Rips complex on a set of $n$ points is $n-1$, and thus all persistence barcodes for $i\geq n-2$ are empty. Similarly as for the enriched nerve, one important property used to derive Equation \ref{E:rips mag} is that the barcodes of the rescaled space $tX$ can be obtained from the rescaled barcodes of $X$. While this is not true for arbitrary simplicial complexes on a metric space $X$,  it holds for simplicial complexes such as Vietoris-Rips, \v{C}ech and alpha.

Unfortunately, as becomes evident from Equation \ref{E:rips mag}, the computation of the Rips magnitude poses some practical issues: even for point clouds with less than $50$ points, it can't be computed with common computational resources. 
One may thus investigate if the persistent magnitude computed using other abstract simplicial complexes, which are computationally more treatable, such as the alpha complex (see Definition \ref{D:alpha}), are still useful to approximate any known fractal dimensions. %
Motivated by these computational observations, in \cite{o2023alpha}  persistent magnitude for the alpha complex of subsets $X$ of Euclidean spaces is introduced, with the name \emph{alpha magnitude}, denoted by $\mathrm{Mag}_\alpha(X)$.
This gives a notion of persistent magnitude that is computable in practice and has opened the door for practical applications.
Furthermore, in \cite{o2023alpha} the definition of alpha magnitude is extended to compact subsets of Euclidean space, and the growth rate of the alpha magnitude function  is studied, leading to the following definition:

\begin{definition}\label{D:amd}
Let $X\subset \mathbb{R}^D$ be a compact subset. 
The \define{alpha magnitude dimension} of $X$ is defined as the number  
\[
\lim_{t \to \infty } \frac{\log(
\mathrm{Mag}_\alpha(tX))}{\log(t)}
\]
whenever this limit exists. 
\end{definition}

For spaces with known fractal structure, such as the Cantor set and the Feigenbaum attractor, alpha magnitude dimension coincides or gives a suitable approximation of the box-counting dimension, respectively.
Such theoretical and heuristic observations led to the formulation of a  conjecture  about the equivalence between the alpha magnitude dimension and the box-counting dimension for compact subsets of Euclidean space \cite{o2023alpha}.

\subsection{Outlook for networks}
We can thus ask: how can we leverage the existing approaches to capture self-similarity of networks? 
We believe that some key properties of topological approaches that make them particularly suitable to estimating fractal dimensions of networks are the following:

\begin{enumerate}
\item  \emph{Ability to capture higher-order information and invariants.} By using simplicial complexes, we can encode higher-order relationships between nodes, thus adding important information. It is worth noting that higher-order relationships can also be captured by hypergraphs, however, simplicial complexes, as topological objects, can naturally  be leveraged by methods such as PH that are able to capture local geometric information, crucial for the estimation of fractal dimensions.  In particular, we would expect  suitable generalisations of PH-dimensions in different degrees to give more fine-grained information about  fractal-like structure present. For instance, we can imagine that different networks with same value for $PH_0$ might take on different values for higher-degree PH-dimensions, which could capture the repeating of higher-order patterns between nodes. 
\item  \emph{Strong theoretical foundations.} PH and magnitude-based methods rest on strong theoretical foundations,  are able to capture intrinsic properties of data, while at the same time being easier to generalise and transfer to different contexts.
\item  \emph{Potential computational advantages.} As mentioned in the previous part of this section, it has been shown by \cite{Jaquette2019FractalDE}, that there are strong computational advantages in using certain types of topological notions of fractal dimensions for point samples and subsets of Euclidean space, e.g., in an extensive benchmarking, the  $PH_0$ dimension   was  shown to be easier to compute, more efficient, while outperforming the box-counting and correlation dimensions. Finally, we add a cautionary remark that for other types of topologically-based notions of dimension, such as the magnitude dimension, or PH-dimensions in degrees larger than $0$, computational challenges need to be further addressed, and we discuss this further in Challenge/Open Problem 4 below.
\end{enumerate}

\section{Challenges and open problems}
As we have seen, methods from   topology and related fields can help in estimating fractal dimensions for subsets of Euclidean space.  
Naively, the existing methods can be directly applied to networks \emph{provided}  one forgets the adjacency information of nodes (i.e., which nodes are connected by an edge) and only retains information about closeness of nodes, encoded in the shortest-path metric structure.

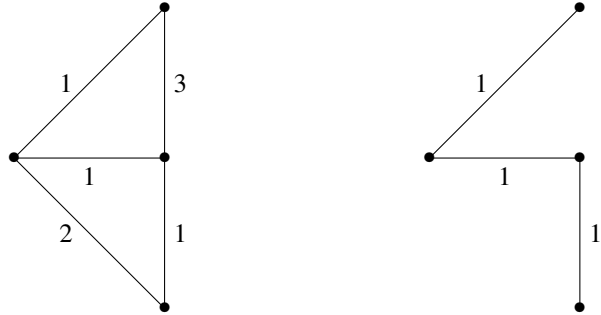
\begin{figure}[h!]
\begin{center}
\begin{tikzpicture}
\node at (0,0) {$\bullet$};
\node at (0,2) {$\bullet$};
\node at (-2,0) {$\bullet$};
\node at (0,-2) {$\bullet$};
\path[-]
(0,0) edge node[right] {$3$} (0,2)
(0,0) edge node[right]{$1$} (0,-2)
(0,-2) edge node[left] {$2\;$} (-2,0)
(-2,0) edge node[below] {$1$} (0,0)
(-2,0) edge node[left] {$1\;$} (0,2);
\end{tikzpicture}
\qquad\qquad\qquad\qquad
\begin{tikzpicture}
\node at (0,0) {$\bullet$};
\node at (0,2) {$\bullet$};
\node at (-2,0) {$\bullet$};
\node at (0,-2) {$\bullet$};
\path[-]
(0,0) edge node[right]{$1$} (0,-2)
(-2,0) edge node[below] {$1$} (0,0)
(-2,0) edge node[left] {$1\;$} (0,2);
\end{tikzpicture}
\end{center}
\caption{Example of two networks with same shortest-path metric structure, where we label the edges with the corresponding weight. We note that the network on the left contains two cliques of $3$ nodes each, while the network on the right is a line network. The edge of length 3 will not be the shortest path between the corresponding edges for the network on the left, even though it is the most direct, because we are considering the shortest \emph{weighted} path, which is of length 2. }\label{F:ex sp metric space}
\end{figure}

Thus a natural question to ask is what information is  lost by the map sending a weighted network $N=(V,E,w)$ to the corresponding  shortest-path metric space $(V,\SP)$. This question has been explored in \cite{HY65}, in which the authors study  properties of networks being uniquely determined by their associated shortest-path metric structure. In particular, no information is lost if the weighted network is a tree. However, the same is not true for more general networks. We give an example of two weighted networks with same shortest-path metric  in Figure \ref{F:ex sp metric space}. In the following we thus make the distinction between methods that would recognise the difference between two networks such as those in  Figure \ref{F:ex sp metric space}, for which we say that they capture information about ``connectedness'' between nodes, and methods that would instead not be able to tell such two networks apart, for which we say that they forget information  about ``connectedness''.

\paragraph{Challenge 1: PH dimensions forget connectedness information between nodes.}
As discussed in Section \ref{SS:PH dim}, PH-based notions of fractal dimensions have strong theoretical foundations, are easy to use, efficient and approximate well the known fractal dimension of subsets of Euclidean space. However, the suggested dimensions in the literature have been defined for measures on bounded metric spaces, and can thus be applied to network provided one forgets the connectedness information between nodes. This poses a problem when applying such methods to capture intrinsic properties of networks,  
 since, as we discussed at the beginning of this section, important network structure is forgotten when passing to its shortest-path metric structure.%

\paragraph{Open Problem 1: Develop  PH dimensions for networks.}
To be able to transfer PH dimensions to networks, without forgetting the connectivity structure, one needs (i) suitable filtrations of simplicial complexes and (ii) a suitable notion of probability measure on a weighted network.

For  point (i), we note that in the simplest case of assigning a simplicial complex to a network, one can build a clique complex whose network backbone is the original network. There are several types of filtrations of this kind that have already been considered in the TDA literature, most of which filter on the edge-weights, and hence are suitable for weighted networks. 
One such filtration is given by so-called \emph{weight-rank clique complexes} \cite{10.1371/journal.pone.0066506}: given a parameter value $\epsilon$ one considers the subnetwork in which one only includes edges with weights bounded above by  $\epsilon$, and one then associates a clique complex to the subnetwork. Repeating this procedure for all possible values of the parameter $\epsilon$ one obtains a filtered simplicial complex.
We refer the reader to \cite{PH-network} for an overview of different types of filtered simplicial complexes that one can associate to a network. 

Finally, for data-driven systems, simplicial complexes are natural models for capturing higher-order interactions, i.e., simultaneous interactions between more than two nodes. For example, in brain networks where neurons are nodes, a $n-$simplex could indicate the simultaneous firing of $(n+1)$ neurons.

For point (ii) we note that studying subsampling strategies for networks is an active area of research in  network science and adjacent fields, see for instance \cite{10.1093/bioinformatics/bth163,4015377}.

\paragraph{Challenge 2: Magnitude forgets connectedness information between nodes.}
Magnitude, as  defined  in Section \ref{SS: magnitude}, is an invariant of a metric space, and magnitude of shortest-path metric spaces associated to (unweighted) networks has been studied in \cite{LEINSTER_2019}. 
Similarly as for PH dimensions, magnitude dimension  forgets important information about the structure of a network.

\paragraph{Open Problem 2: Develop a notion of magnitude that takes into account adjacency information between nodes in a network.}
The definition of magnitude of a metric space arises as a special instance of a more general notion of magnitude of a category developed by Tom Leinster \cite{TL08}, of which the cardinality of a set, and the Euler characteristic of a topological space can be seen as specific instantiations. Similarly as the Euler characteristic of a topological space can be recovered by computing alternating sums of the ranks of its singular homology groups, magnitude of a finite metric space is related to magnitude homology \cite{HW14,LS21}. In turn, magnitude homology of a finite metric space has been shown to be equivalent to the persistent homology of  a specific type of filtered simplicial set, called \emph{the enriched nerve} of the metric space \cite{otter2021magnitude}. Thus a natural question that arises is: can one define a network-analogue of the enriched nerve of a metric space that preserves connectedness relationships between nodes?

\paragraph{Challenge 3: Existing persistent magnitude functions cannot be used directly  for networks.}
There are currently two notions of persistent magnitude functions that have been studied in the literature, and from which one can derive appropriate notions of dimension\footnote{We note that of these, only alpha persistent magnitude has been used to define a corresponding dimension  (see Definition \ref{D:amd}), however one could in an analogous way also define Rips persistent magnitude dimension.}: Rips persistent magnitude  and alpha persistent magnitude. The first can be defined for any metric space, while for the second we need for the metric space to be embedded into Euclidean space. Rips persistent magnitude  can in theory be computed for any metric space, and thus also for the shortest-path metric space associated to a network; however, in practice this quantity is not computable for networks having more than a few nodes, and thus cannot be used to estimate fractal dimensions of networks. Most networks cannot be embedded into Euclidean space, thus alpha magnitude dimension, while better suited to computations than Rips magnitude dimension, does not provide a suitable method to estimate dimensions of networks.  \\

\paragraph{Open Problem 3: Develop notions of persistent magnitude dimensions that preserve the network structure.}
Since the definition of persistent magnitude dimension depends on a specific choice of filtration, a solution to Challenge 4 is given by associating with a given network a simplicial complex that preserves the connectedness relationship between the nodes, while at the same time also capturing  the information given by the edge weights, by using one of the filtrations of simplicial complexes discussed in Open Problem 1.

\paragraph{Challenge 4: PH-dimensions and magnitude  can be computationally expensive.}
While as we discussed above, $PH_0$-dimension has been proven to be computationally efficient, the same is not true for $PH_i$ for $i\geq 1$. 
One way to understand the computational complexity of PH-dimensions is to consider what the computational complexity of the PH computation is over the different samples needed to compute empirical expectations in the definition of PH-dimension. 
If $N$ is the number of simplices in a given filtered simplicial complex, then the computational complexity of the PH computation can be understood as being $\mathcal{O}(N^{2.8})$, although in practice it is often much lower. We refer the reader to \cite{otter2017roadmap} (and references therein)
for more details on computational complexity of PH algorithms.

For what concerns magnitude, 
the standard method of computing magnitude for a set of $n$ points requires inverting an $n \times n$ matrix. The best known lower bound for matrix multiplication and inversion is $\Omega(n^2\log n)$ \cite{raz2002complexity}. 
Hence, since the computation of magnitude needs all pairs of input points, it becomes expensive for large datasets.

Some recent efforts have been made to reduce the computational complexity of magnitude by framing the magnitude computation in the language of convex optimisation \cite{andreeva2024approximating}. Another proposed approach was to use carefully selected subsamples\cite{andreeva2024approximating}. Empirical results demonstrate that these approaches work well in practice for a number of synthetic and real datasets. However, due to the fact that magnitude is not submodular \cite{andreeva2024approximating}, there are no theoretical guarantees in this direction. This leads us to the following open problem.

\paragraph{Open Problem 4: Develop efficient algorithms for the computation of PH-dimensions or magnitude.}
As discussed in Challenge 4, the computational cost of magnitude and persistent homology in high degrees can make it prohibitive for magnitude or PH-dimensions to estimate fractal dimensions for many spaces or networks of interest.

For what concerns PH-dimensions, we note that 
filtered simplicial complexes associated to networks have the advantage of being potentially sparser (i.e., having fewer simplices) than simplicial complexes associated to metric spaces. For instance, the Vietoris-Rips complex on $n$ points in a metric space has $2^n$ simplices. If one is only interested in computing PH in degree $i$, then one may construct the Vietoris-Rips complex only up to simplices of dimension $i+1$, and thus the resulting complex will contain approximately $n^{i+1}$ many simplices. On the other hand, in a filtered simplicial complex such as the weight-rank clique filtration, the number of simplices in a given dimension $i$ is determined by the number of $i$-cliques present in the starting network.
Thus, while PH-dimensions for finite point clouds might be difficult to compute in degrees higher than $0$, the same might not be true for a network counterpart of this notion of dimension. This makes this a very interesting direction to be further explored. 

For what concerns magnitude, 
 there have been recent advances trying to approximate the magnitude of a space by finding appropriate  sparsification methods yielding  subsamples small in size, with magnitude close to that of the original space \cite{andreeva2024approximating}. However, this approach does not provide theoretical guarantees, thus further effort is needed in this direction.

We also note that there is a notion closely linked to magnitude, called \emph{spread} \cite{willerton2015}, which has been shown to be more tractable, and for which there is heuristic evidence of it being able to estimate box-counting dimension of some fractals such as the Cantor Set or Koch curve. Thus, further investigating this notion might present another avenue to overcome the difficulty in computing magnitude.

\paragraph{Further open problem}

We conclude this paper by discussing one further open problem of a wider scope that, while not being strictly related to networks, would also need to be addressed as part of the topological study of fractal dimensions for networks.

We think that there is a need for  a systematic comparison of the various existing topological fractal dimensions, against existing methods such as box-counting or correlation dimensions. To our knowledge, the only benchmarking relating topological fractal dimensions and more classical estimations of fractal dimensions has been performed for the PH dimension.
We believe that a more comprehensive benchmarking, which would include magnitude dimension and persistent magnitude dimension for point clouds would help elucidate the advantages of the different approaches, and would pave the way for a similar comparison for fractal dimensions for networks. Furthermore, conducting a more thorough investigation of the spread dimension, which seems to have been under-explored both theoretically and empirically, will help provide more clarity on its properties, and could be a useful addition to the benchmarking suite of fractal dimensions.
Finally, another notion relating to magnitude are fractional curvature measures, which generalise fractal dimensions to intrinsic volumes, such as the Euler characteristic, see \cite{Winter}. It would thus be of interest to  investigate how this notion compares to magnitude or spread dimension.

\section*{Acknowledgements}
We thank the organisers of the 3rd Women in Computational Topology (WinCompTop), who brought our group together in this collaboration, as well as the AWM for the generous financial support during the WinCompTop week, and the Bernoulli Center for their hospitality.  RA, NO and ET thank the ICMS as well as the Maxwell Foundation for their support and hospitality during their ICMS Research in Group stay in Edinburgh in July 2024. HCP was supported by the UNAM Posdoctoral Program (POSDOC) and DGAPA-PAPIIT IN114323.
We thank the anonymous Reviewers for their constructive feedback that has helped improve our manuscript.

\bibliographystyle{alpha}
 \bibliography{references}
\end{document}